\documentclass[12pt]{article}
\usepackage[colorlinks,
            linkcolor=blue,
            anchorcolor=green,
            citecolor=blue
            ]{hyperref}
\usepackage{mathrsfs,amssymb,amsfonts}
\usepackage{amsmath,amssymb,latexsym,color}
\usepackage[mathscr]{eucal}
\usepackage{algorithm,algorithmic}
\usepackage[numbers,sort&compress]{natbib}
\usepackage{graphicx}
\usepackage{comment}
\usepackage{lipsum}
\def\thefootnote{\fnsymbol{footnote}}
\newtheorem{thm}{Theorem}[section]

\newtheorem{lemma}[thm]{Lemma}
\newtheorem{cor}[thm]{Corollary}
\newtheorem{example}[thm]{Example}

\newtheorem{problem}[thm]{Problem}

\newtheorem{remark}[thm]{Remark}
\newtheorem{Notation}[thm]{Notation}

\newcommand{\proof}{{\it Proof.\quad}}
\newcommand{\qed}{\hfill\Box\medskip}
\renewcommand{\thefootnote}{\arabic{footnote}}

\newcommand{\diam}{\mathop{\mathrm{diam}}}

\newcommand{\Sub}{\mathop{\mathrm{Sub}}}

\usepackage{authblk}

\begin{document}
\renewcommand{\abovewithdelims}[2]{
\genfrac{[}{]}{0pt}{}{#1}{#2}}

\title{\bf The diameter and dominating sets of the difference graph of a nilpotent group
}

\author[1]{Xuanlong Ma}
\author[2]{Samir Zahirovi\' c}
\author[3]{Katarina \v Zigerovi\' c}
\affil[1]{\small School of Science, Xi'an Shiyou University, Xi'an 710065, P.R. China}
\affil[2]{\small University of Novi Sad, Faculty of Sciences, Department of Mathematics and Informatics, Serbia}
\affil[3]{\small Mathematical Institute of the Serbian Academy of Sciences and Arts, Serbia}

 \date{}
 \maketitle
\newcommand\blfootnote[1]{%
\begingroup
\renewcommand\thefootnote{}\footnote{#1}%
\addtocounter{footnote}{-1}%
\endgroup
}

\begin{abstract}
Given a finite group $G$, the difference graph of $G$, denoted by $\mathcal{D}(G)$, is the difference of the enhanced power graph of $G$ and the power graph of $G$, with all isolated vertices removed.
This paper mainly studies the dominating sets of the difference graph of a finite group. In particular, we prove that the diameter of the difference graph of a nilpotent group has an upper bound of $4$. Furthermore, we generalize and refine the result by Biswas et al. \cite[Theorem 5.3]{BC22} by
classifying all nilpotent groups whose difference graph has diameter $k$, for each $k\le 4$.


\end{abstract}

\blfootnote{{\em E-mail addresses:}\\
\url{xuanlma@xsyu.edu.cn} (X. Ma),\\
\url{samir@dmi.uns.ac.rs} (S. Zahirovi\' c),\\
\url{katarina.zigerovic@mi.sanu.ac.rs} (K. \v Zigerovi\' c)
}

\section{Introduction}

For an algebraic structure (group or ring), one can define
a graph on this algebraic structure and study its graph representations,  which is an interesting research topic in algebraic graph theory and has valuable applications.
Such as, the Cayley graphs of groups as classifiers for data mining \cite{KeR}.

One of the well-known graph representations on groups is the
power graphs. Given a group $G$,
the {\em power graph} of $G$, denoted by $\mathcal{P}(G)$, is the simple graph whose vertex set is $G$, such that two vertices are adjacent if one is a power of the other.
Kelarev and Quinn \cite{KQ00} first introduced the directed power graph of a group. In 2009, Chakrabarty {\em et al.} \cite{CGS} introduced
the undirected power graph of a group, which is also called the power graph of a group.
In past two decades, the study of power graphs of groups has been growing (see the two survey papers \cite{AKC,KSC}).
The {\em enhanced power graph} of a group $G$, denoted by $\mathcal{P}_E(G)$, is the simple graph with vertex set $G$, where two distinct vertices are adjacent if they generate a cyclic subgroup.
Aalipour {\em et al.} \cite{Aa17} first introduced the term enhanced power graph.
In \cite{Aa17}, the authors studied the enhanced power graph in order to measure how close the undirected power graph is to the commuting graph.
A detailed list of results and open problems related to enhanced power graph of a group can be found in \cite{MaW}.

In 2022, Biswas, Cameron {\em et al.} \cite{BC22} first introduced the difference graph of a finite group. In fact, it was observed earlier in \cite{cameron-graphs-defined-on-groups} by Cameron that the power graph, the enhanced power graph, the commuting graph, and several other graphs associated with the groups form a hierarchy, in the sense that each of them is a subgraph of the next. He suggested investigating properties of differences of these graphs, and the difference of the enhanced power graph and the power graph is the one that has been studied the most. Given a finite group $G$, the {\em difference graph} of $G$, denoted by $\mathcal{D}(G)$,
is the difference of the
enhanced power graph of $G$ and the power graph of $G$, where all isolated vertices are removed. That is, $\mathcal{D}(G)$ is the graph obtained by removing all isolated vertices in  $\mathcal{P}_E(G)-\mathcal{P}(G)$. Recently, Kumar {\em et al.} \cite{Ku} characterized the finite groups whose difference graph is a chordal graph, a star graph, a dominatable graph, a threshold graph, and a split graph. Bo\v snjak, Madar\' asz and Zahirovi\' c \cite{bosnjak-madarasz-zahirovic} proved that, if two abelian groups of orders divisible by only two primes have isomorphic difference graphs, then they are isomorphic.

In Section \ref{domination-one}, we gave a description of finite groups whose domination number is $1$. In Section \ref{cyclic groups}, we turn our attention to cyclic groups; we provide a full description of dominating sets of minimal size of the difference graph of a cyclic group. Thus, we determine classes of finite cyclic groups $\mathbb Z_n$ such that $\gamma(\mathcal D(\mathbb Z_n))=1$ and $\gamma(\mathcal D(\mathbb Z_n))=2$, as well as classes of finite cyclic groups such that $\diam(\mathcal D(\mathbb Z_n))=2$ and $\diam(\mathcal D(\mathbb Z_n))=3$. In Section \ref{non-cyclic nilpotent}, we prove that the diameter of the difference graph fo a finite nilpotent group is at most $4$. We construct some domination sets of the difference graph of a finite non-cyclic nilpotent group. Later in the section, we use that to describe nilpotent groups whose difference graph has diameter $2$, $3$ and $4$.

\section{Preliminaries}

This section introduces some fundamental definitions from group theory and graph theory that will be used throughout the paper.

All groups considered are assumed to be finite. We denote groups by capital letters, usually  $G$, and its identity element by $e$. For an element $g\in G$, the {\em order} of $g$ is denoted by $o(g)$. In particular, if $o(g)=2$, then $g$ is called an {\em involution}.
The {\em exponent} of $G$, denoted by $\exp(G)$, is the least common multiple of the orders of all elements in $G$. Equivalently, it is the smallest positive integer $t$ such that $g^t=e$, for every $g\in G$.

A {\em maximal cyclic subgroup} of $G$ is a cyclic subgroup which is not a proper subgroup of some other cyclic subgroup of $G$.
The set of all maximal cyclic
subgroups of $G$ is denoted by $\mathcal{M}(G)$.
Clearly, $G$ is cyclic if and only if
$\mathcal{M}(G)=\{G\}$, or equivalently, if and only if $|\mathcal{M}(G)|=1$.

The cyclic group of order $n$ is denoted by $\mathbb{Z}_n$.
Given a prime $p$, the elementary abelian $p$-group of order $p^n$ is the direct product of $n$ copies of $\mathbb{Z}_p$, and is denoted by $\mathbb{Z}_p^n$.
It is well known that a finite group $G$ is {\em nilpotent} if and only if every Sylow $p$-subgroup of $G$ is normal, which is also equivalent to saying that $G$ is the direct product of its Sylow subgroups. In particular, if $G$ is nilpotent and $g,h\in G$ with $\gcd(o(g),o(h))=1$, then $gh=hg$.

All graphs considered in this paper are assumed to be simple, that is, undirected graphs without loops or multiple edges.
For a graph $\Gamma$, we denote its vertex set by $V(\Gamma)$, and its edge set by $E(\Gamma)$. If two distinct vertices $x$ and $y$ of $\Gamma$ are adjacent, we denote this by $x\sim y$. The {\em distance} between $x$ and $y$ in $\Gamma$ is the length of a shortest path from $x$ to $y$, and is denoted by $d(x,y)$. The {\em diameter} of $\Gamma$, denoted by $\diam (\Gamma)$, is the greatest distance  between any two vertices in $\Gamma$. A subset $D\subseteq V(\Gamma)$ is called a {\em dominating set} of a graph $\Gamma$ if every vertex $v\in V(\Gamma)$ is either in $D$ or is adjacent to at least one vertex in $D$. The minimum size of dominating sets  of  $\Gamma$, denoted by $\gamma (\Gamma)$, is called the {\em domination number} of $\Gamma$.

\section{Difference graphs with domination number one}\label{domination-one}

Given a finite group $G$,
the difference graph of $G$, denoted by $\mathcal{D}(G)$,
is defined as the difference of the
enhanced power graph of $G$ and the power graph of $G$, where all isolated vertices are removed. Namely, $\mathcal{D}(G)$ is the graph obtained by removing all isolated vertices in $\mathcal{P}_E(G)-\mathcal{P}(G)$.

\begin{lemma}{\rm (\cite[Proposition 4.2]{BC22})}\label{two-con}
Let $a,b\in G$ be two distinct non-trivial elements. If $a$ and $b$ commute and $\gcd(o(a),o(b))=1$, then $a\sim b$ in $\mathcal{D}(G)$.
\end{lemma}

It is clear that the identity element $e\in G$ is always isolated vertex of $\mathcal{P}_E(G)-\mathcal{P}(G)$, and therefore $e\notin V(\mathcal{D}(G))$.
The following lemma characterizes all non-trivial isolated vertices of $\mathcal{P}_E(G)-\mathcal{P}(G)$.

\begin{lemma}{\rm (\cite[Proposition 2.1]{BC22})}\label{iso-vertex}
Let $G$ be a group of order at least $2$.
A non-trivial element $g\in G$ is isolated in $\mathcal{P}_E(G)-\mathcal{P}(G)$ if and only if either
$\langle g\rangle$ is a maximal cyclic subgroup of $G$, or every cyclic subgroup of $G$ containing $g$ has prime-power order.
\end{lemma}

We now introduce a special class of groups that play a key role in the structure of difference graphs with domination number one.

Let $\Phi$ denote the class of finite groups $G$ satisfying the following conditions:
\begin{itemize}
\item[{\rm (a)}] There exists a unique involution in  $G$ that commutes with at least one element of odd prime order. Let us denote this involution by $u$;

\item[{\rm (b)}] For every non-trivial element $g\in G$ of odd order, one of the following holds: $\langle g\rangle$ is a maximal cyclic subgroup of $G$; every cyclic subgroup of $G$ containing $g$ has prime-power order; $g$ has prime order and commutes with $u$;

\item[{\rm (c)}] For every element $g\in G\setminus\{u\}$ of even order,  either $\langle g\rangle$ is a maximal cyclic subgroup of $G$, or every cyclic subgroup of $G$ containing $g$ has prime-power order.
\end{itemize}
A group $G$ is called a {\em $\Phi$-group} if $G\in \Phi$. Note that $\Phi$ is not empty; for example, $\mathbb{Z}_2\times \mathbb{Z}_q^n\in \Phi$ for any odd prime $q$ and any positive integer $n$.

\begin{thm}\label{d-n-1}
$\gamma(\mathcal{D}(G))=1$ if and only if $G$ is a $\Phi$-group. Moreover, if $\gamma(\mathcal{D}(G))=1$, then $\{u\}$ is the unique dominating set of size one in $\mathcal{D}(G)$, where $u$ is the unique involution in $G$ such that there exists at least one element of odd prime order that commutes with $u$.
\end{thm}
\proof
Assume first that $\gamma(\mathcal{D}(G))=1$, and let $\{u\}$ be a dominating set of $\mathcal{D}(G)$. If $u$ is not an involution, then $u^{-1}\ne u$ and $u^{-1}\in V(\mathcal{D}(G))$, so $u^{-1}$
is not adjacent to $u$ in $\mathcal{D}(G)$, which is impossible.
As a result, we have that $u$ is an involution.
Since $u\in V(\mathcal{D}(G))$, it follows from Lemma~\ref{iso-vertex} that there exists a cyclic subgroup of $G$ containing $u$ whose order is not a prime power. This implies that $u$ commutes with an element of odd prime order.

Let $g\in G\setminus\{u\}$ be of even order. If $g\sim u$ in $\mathcal{P}_E(G)$, then $g\sim u$ also in $\mathcal{P}(G)$. It follows that $g$ is not adjacent to $u$ in $\mathcal{P}_E(G)-\mathcal{P}(G)$,  so $g$ is an isolated vertex of $\mathcal{P}_E(G)-\mathcal{P}(G)$. By Lemma~\ref{iso-vertex}, condition (c) is satisfied.
This also means that (a) in the definition of a $\Phi$-group holds.

Now, let $h\in G$ be a non-trivial element of odd order.
If $h$ is an isolated vertex of $\mathcal{P}_E(G)-\mathcal{P}(G)$, then by Lemma~\ref{iso-vertex}, either
$\langle h\rangle$ is a maximal cyclic subgroup of $G$, or every cyclic subgroup of $G$ containing $h$ has prime-power order. Thus, in the following, we may assume that $h\in V(\mathcal{D}(G))$. Then we have that $h\sim u$ in $\mathcal{D}(G)$, so $\langle h,u\rangle$ is cyclic, and $h$ and $u$ commute. Suppose, for contradiction, that $o(h)$ is not prime. Then there exists an element $c\in \langle h\rangle$ of prime order, so $c\sim u$, and $cu$ has even order. Consequently, Lemma~\ref{iso-vertex} implies that $cu\in V(\mathcal{D}(G))$, so $cu\sim u$ in $\mathcal{D}(G)$. However, this is impossible, as $cu\sim u$ in $\mathcal{P}(G)$. Hence, $o(h)$ must be the prime, and (b) is satisfied. Therefore, $G\in \Phi$.

Conversely, suppose that $G\in \Phi$. Let $u$ be the unique involution in $G$ such that there exists an element of odd prime order that commutes with $u$. Then, by Lemma~\ref{two-con}, $u\in V(\mathcal{D}(G))$.
Take any $v\in V(\mathcal{D}(G))\setminus\{u\}$. By (c) of the definition of a $\Phi$-group and Lemma~\ref{iso-vertex}, $o(v)$ must be odd. Note that $v\ne e$. Then, by (b), we conclude that $v$ has a prime order and commutes with $u$, so $v\sim u$ in $\mathcal{D}(G)$.
Hence, $\{u\}$ is a dominating set of $\mathcal{D}(G)$ and $\gamma(\mathcal{D}(G))=1$, as desired.
$\qed$

The next result classifies all finite nilpotent groups whose
difference graphs have domination number one.

\begin{cor}\label{nil-1}
Let $G$ be a nilpotent group. Then $\gamma(\mathcal{D}(G))=1$ if and only if
$$G\cong \mathbb{Z}_2\times Q,$$ where $Q$ is a group
with exponent $q$, for some odd prime $q$.
\end{cor}
\proof
Clearly, for an odd prime $q$, if $Q$ is a group with $\exp(Q)=q$, then $\mathbb{Z}_2\times Q$ is a nilpotent $\Phi$-group, and so $\gamma(\mathcal{D}(G))=1$ by Theorem~\ref{d-n-1}.

Now, let $G$ be a nilpotent $\Phi$-group with a unique involution $u$. Then $|G|$ must be divisible by $2$. Suppose that $|G|$ has at least three distinct prime divisors. Then there exists an element $x\in G$ of order $pq$, where $p$ and $q$ are distinct odd prime divisors of $|G|$. But now, $ux$ has order $2pq$ and $x\in \langle ux\rangle$,
which is a contradiction by (b) in the definition of a $\Phi$-group.

We conclude that $|G|$ has at most two prime divisors. Thus, $G\cong P\times Q$, where $P$ is a $2$-group and $Q$ is a $q$-group for some odd prime $q$.
Moreover, by (a) of the definition of a $\Phi$-group again, we see that $P$ has a unique involution, say $u$. We claim that $P$ has no elements of order $4$. In fact, if $P$ has an element $a$ of order $4$, then $G$ has a cyclic subgroup of order $4q$ containing $a$, which is impossible by (c). Therefore, $P\cong \mathbb{Z}_2$.
To complete the proof, note that if $Q$ contains an element $x$ of order $q^2$, then $\langle ux\rangle$ has order $2q^2$ and contains $x$, again contradicting (b). Hence,  $\exp(Q)=q$, as required.
$\qed$

\begin{cor}\label{abe-1}
Let $G$ be an abelian group. Then $\gamma(\mathcal{D}(G))=1$ if and only if
$$G\cong \mathbb{Z}_2\times \mathbb{Z}_q^n,$$
where $q$ is an odd prime and $n$ is a positive integer.
\end{cor}

\section{The diameter and dominating sets of the difference graph of a cyclic group}\label{cyclic groups}

In \cite{Aa17}, Aalipour {\em et al}. showed that the enhanced power graph and the power graph of a finite group $G$ are equal if and only if every cyclic subgroup of $G$ has prime power order. Thus, if $G$ is a $p$-group, then the difference graph $\mathcal{D}(G)$ is not defined. For this reason, we restrict our attention to the cyclic groups $\mathbb{Z}_n$, where $n$ is not a prime power. Note that, in this case, we have that
$V(\mathcal{D}(\mathbb{Z}_n))=\{g\in \mathbb{Z}_n\setminus\{e\}: o(g)\ne n\}$, and
two distinct vertices $a,b\in V(\mathcal{D}(\mathbb{Z}_n))$ are adjacent if and only if $\gcd(o(a),o(b))\notin \{o(a),o(b)\}$.

 \begin{lemma}\label{cyc-domn}
Let $n$ be a positive integer which is not a prime power. Then
\begin{equation}\label{cyc-dn}
\gamma(\mathcal{D}(\mathbb{Z}_n))=\left\{
                          \begin{array}{ll}
                          1, & \hbox{if $n=2q$ for some odd prime $q$;} \\
                          2, & \hbox{otherwise.}
                          \end{array}
                        \right.
\end{equation}
Furthermore, if $n=2q$ for an odd prime $q$, then the unique dominating set of size $1$ is $\{u\}$, where $u$ denotes the unique involution of $\mathbb Z_n$. Otherwise, $\{a,b\}$ is a dominating set of size $2$ if and only if $\gcd(o(a),o(b))=1$ and $o(a)o(b)=n$, $a$ and $b$ are adjacent and $\{a,b\}$ is a dominating set in $\mathcal{D}(\mathbb{Z}_n)$.
\end{lemma}
\proof
 If $n=2q$ for some odd prime $q$, it is easy to see that $\mathcal{D}(\mathbb{Z}_n)\cong K_{1,q-1}$, where the unique involution of $\mathbb Z$ is the unique vertex of degree $n-1$ in $\mathcal D(\mathbb Z_n)$.

 Now, suppose that $n$ is not a product of $2$ and an odd prime. Let $a,b\in \mathbb{Z}_n$ be such that $\gcd(o(a),o(b))=1$ and $o(a)o(b)=n$. Then Corollary~\ref{abe-1} also implies that $\gamma(\mathcal{D}(\mathbb{Z}_n))\ge 2$.
 Note that
$$N(a)=\{a'\in \mathbb{Z}_n\setminus \langle a\rangle: a'\neq e, o(a)\nmid o(a')\},$$
$$N(b)=\{b'\in \mathbb{Z}_n\setminus \langle b\rangle : b'\neq e, o(b)\nmid o(b')\}.$$

 Consider any $x\in V(\mathcal{D}(\mathbb{Z}_n))\setminus\{a,b\}$.
It is easy to see that if, without loss of generality, $x\in \langle a\rangle$, then $x\sim b$. Now suppose that $x\in \mathbb{Z}_n \setminus (\langle a\rangle \cup \langle b\rangle)$. Since $o(x)\neq n=o(a)o(b)$, either $o(a)\nmid o(x)$ or $o(b)\nmid o(x)$. If, without loss of generality, $o(a)\nmid o(x)$, then $a\sim x$. Therefore, $\{a,b\}$ is a dominating set of $\mathcal{D}(\mathbb{Z}_n)$, and therefore $\gamma(\mathcal{D}(\mathbb{Z}_n))=2$, as desired.

It remains for us to show that, whenever $\{a,b\}$ is a dominating set of $\mathcal D(\mathbb Z_n)$, then $o(a)o(b)=n$ and $a$ and $b$ have relatively prime orders. Let $k=\gcd(o(a),o(b))$. If $k>1$, then none of $a$ and $b$ is adjacent to any of the elements of order $k$ of $\mathbb Z_n$. Therefore, $o(a)$ and $o(b)$ are relatively prime. Finally, if $o(a)o(b)<n$, then $\mathcal D(\mathbb Z_n)$ contains elements of order $o(a)o(b)$ of $\mathbb Z_n$, and neither $a$ nor $b$ is adjacent to an element of order $o(a)o(b)$. Thus, the lemma has been proved.
$\qed$

\begin{thm}\label{cyc-diam}
Let $n$ be a positive integer which is not a prime power. Then
\begin{equation}\label{cyc-dn}
\diam(\mathcal{D}(\mathbb{Z}_n))=\left\{
                          \begin{array}{ll}
                          2, & \hbox{if $n$ is a square-free number;} \\
                          3, & \hbox{otherwise.}
                          \end{array}
                        \right.
\end{equation}
In particular, $\mathcal{D}(\mathbb{Z}_n)$ is connected.
\end{thm}
\proof Observing that $\mathcal{D}(\mathbb{Z}_n)$ is not a complete graph, Lemma~\ref{cyc-domn} implies that $2\leq \diam(\mathcal{D}(\mathbb{Z}_n))\leq 3$.

Suppose that $n$ is a square-free number, and let
$x,y\in V(\mathcal{D}(\mathbb{Z}_n))$ where $x\ne y$. If $o(x)\nmid o(y)$ and $o(y)\nmid o(x)$, then we have $x\sim y$. Therefore, without loss of generality, we may assume that $o(x)\mid o(y)$.
Since $o(y)\neq n$ and $n$ is not a prime power, there exists at least one prime divisor $p$ of $n$ such that $p\nmid o(y)$.
Let $z\in \mathbb{Z}_n$ be an element of order $p$. Then $z\in V(\mathcal{D}(\mathbb{Z}_n))$, and  $x\sim z \sim y$. Hence, $d(x,y)=2$, and so in this case we conclude that $\diam(\mathcal{D}(\mathbb{Z}_n))=2$, as desired.

Now suppose that there exists a prime number $p$ such that $p^2 \mid n$. Then, we may write $n=p^mk$, where $m\geq 2$, $p^{m+1}\nmid n$ and $k$ is a positive integer.
Let $a,b\in \mathbb{Z}_n$ be such that
$$
o(a)=p,~~o(b)=p^{m-1}k.
$$
Since $o(a)\mid o(b)$, we have that
$a,b\in V(\mathcal{D}(\mathbb{Z}_n))$ and $a\nsim b$. In the following, we claim that $d(a,b)\ge 3$. Suppose for a contradiction that there exists an element $c\in V(\mathcal{D}(\mathbb{Z}_n))$ such that $a\sim c \sim b$ in $\mathcal{D}(\mathbb{Z}_n)$. Since $o(c)\nmid o(b)$ and $o(b)\nmid o(c)$, it follows that $p^{m}\mid o(c)$. Therefore, $o(a)\mid o(c)$, implying that $a$ and $c$ are not adjacent in $\mathcal{D}(\mathbb{Z}_n)$, which is a contradiction. Thus, the claim holds, i.e., $d(a,b)\ge 3$, which implies that $\diam(\mathcal{D}(\mathbb{Z}_n))\ge 3$. Recalling that $\diam(\mathcal{D}(\mathbb{Z}_n))\leq 3$, we conclude that $\diam(\mathcal{D}(\mathbb{Z}_n))=3$ in this case.
$\qed$

We use $D_{2n}$ to denote the dihedral group of order $2n$, where $n\ge 3$.
In general, $D_{2n}$ has the following presentation:
\begin{equation}\label{d2n}
D_{2n}=\langle a,b: b^2=a^n=(ab)^2=e\rangle.
\end{equation}
Remark that $D_{2n}$ is non-abelian.
Note that $\{\langle a\rangle,\{ab,a^2b,\ldots,a^{n-1}b,b\}\}$ is a partition of $D_{2n}$, and $a^ib$ is an involution for each $0\le i \le n-1$. Hence, for every $0\le i \le n-1$,
$\langle a^ib \rangle\in \mathcal{M}(D_{2n})$.

We now observe that in $\mathcal{D}(D_{2n})$ elements of the form $a^ib$ are isolated vertices, and the remaining vertices are elements of the cyclic group generated by $a$, which is the cyclic group of order $n$. Therefore,
$\mathcal{D}(D_{2n})\cong \mathcal{D}(\mathbb{Z}_n)$, for any $n\ge 3$.

Now, as a corollary of Theorem~\ref{cyc-diam}, we have the following result which extends  \cite[Theorem 5.5]{BC22}.

\begin{cor}
Assume that $D_{2n}$ is the dihedral group as presented in \eqref{d2n}, where $n\ge 3$. If $n$ is a prime power, then $\mathcal{D}(D_{2n})$ is not defined. Otherwise,
\begin{equation*}
\diam(\mathcal{D}(D_{2n}))=\left\{
                          \begin{array}{ll}
                          2, & \hbox{if $n$ is a square-free number;} \\
                          3, & \hbox{otherwise.}
                          \end{array}
                        \right.
\end{equation*}
\end{cor}

For a positive integer $m\ge 2$,
the generalized quaternion group with order $4m$, denoted by $Q_{4m}$, is defined by Johnson (see \cite[pp. 44--45]{Jon}) and
has the following presentation:
\begin{equation}\label{q4m}
Q_{4m}=\langle x,y: x^m=y^2, x^{2m}=y^4=e, y^{-1}xy=x^{-1}\rangle.
\end{equation}
Recall that $Q_{8}$ is the usual quaternion group of order $8$, namely,
$$Q_{8}\cong \langle i,j,k: i^2=j^2=k^2=ijk=-1\rangle.$$
It is easy to check that $Q_{4m}$ has a unique
involution $x^m=y^2$.
Remark that
\begin{equation}\label{q4m-1}
Q_{4m}=\langle x\rangle \cup \{x^iy: 1\le i \le 2m\},~~
o(x^iy)=4 \text{ for any } 1\le i \le 2m,
\end{equation}
and
\begin{equation}\label{q4m-2}
\mathcal{M}(Q_{4m})=
\{\langle x\rangle,\langle xy\rangle,
\langle x^2y\rangle,\ldots,\langle x^my\rangle\},~~~~x^m\in \bigcap_{M\in \mathcal{M}(Q_{4m})}M.
\end{equation}
Thus, by \eqref{q4m-1} and \eqref{q4m-2} it follows that $\mathcal{D}(Q_{4m})\cong \mathcal{D}(\mathbb{Z}_{2m})$, for any $m\ge 2$.

Again, as a corollary of Theorem~\ref{cyc-diam}, we have the following.

\begin{cor}
Assume that $Q_{4m}$ is the generalized quaternion group as presented in \eqref{q4m}, where $m\ge 2$. If $m$ is a power of $2$, then $\mathcal{D}(Q_{4m})$ is not defined. Otherwise,
\begin{equation*}
\diam(\mathcal{D}(Q_{4m}))=\left\{
                          \begin{array}{ll}
                          2, & \hbox{if $2m$ is a square-free number;} \\
                          3, & \hbox{otherwise.}
                          \end{array}
                        \right.
\end{equation*}
\end{cor}

\section{The diameter of the difference graph of a non-cyclic nilpotent group}\label{non-cyclic nilpotent}

In the previous section, we described dominating sets of cyclic groups, and we used that information to determine the diameter and the domination number of the difference graph of a cyclic group. In this section, we continue our study by focusing on the remaining finite nilpotent groups. To that end, in this section, we always consider non-cyclic nilpotent groups.

Let $G$ be a non-cyclic nilpotent group which is not a $p$-group.
Since a finite nilpotent group is the direct product of its Sylow subgroups, we may assume that
\begin{equation}\label{nil-fj}
G=P_1\times P_2 \times \cdots \times P_t \times \mathbb{Z}_{n},
\end{equation}
where each $P_i$ is the non-cyclic Sylow $p_i$-subgroup of $G$, $p_1,\ldots, p_t$ are pairwise distinct primes, and $p_i\nmid n$ for any $1\le i \le t$. Note that in \eqref{nil-fj},
$n$ may be equal to $1$, and in which case $t\ge 2$.

In order to better understand the dominating sets of the difference graph of a non-cyclic nilpotent group, we will start by taking our focus on non-cyclic $p$-groups. That will help us in determining the diameter of the difference graph of non-cyclic nilpotent group.
Let $P$ be a non-cyclic $p$-group, where $p$ is a prime. On the set of all maximal cyclic subgroups of $P$, we define domination relation as follows:
If  $M,N\in \mathcal{M}(P)$ such that $p^2\le|M|\le |N|$ and $|M\cap N|=|M|/p$, then we say that $M$ is dominated by $N$, and we denote this by $M \ll N$.

\begin{lemma}\label{ma-lem-1}
Let $P$ be a non-cyclic $p$-group, where $p$ is a prime. If $M_1,M_2$ and $M_3$ are three distinct maximal cyclic subgroups of $ P$ such that $M_1 \ll M_2$ and $M_2 \ll M_3$, then $M_1 \ll M_3$.
\end{lemma}
\proof
Note that $|M_1|\le |M_2|\le |M_3|$,
$|M_1\cap M_2|=|M_1|/p$ and $|M_2\cap M_3|=|M_2|/p$. Naturally, $|M_2|/p \ge |M_1|/p$. Notice that $M_1\cap M_2$ and $M_2\cap M_3$ are maximal proper subgroups of $M_1$ and $M_2$, respectively. Therefore, since $M_1\cap M_2$ and $M_2\cap M_3$ are both subgroups of cyclic group $M_2$, then $M_1\cap M_2\subseteq M_2\cap M_3$. Thus, $M_1\cap M_2\subseteq M_1\cap M_3$, and because both $M_1\cap M_2$ and $M_1\cap M_3$ are maximal proper subgroups of maximal cyclic subgroup $M_1$, it follows that $M_1\cap M_2=M_1\cap M_3$, which implies that $\lvert M_1\cap M_3\rvert=\lvert M_1\rvert/p$. Thus, $M_1 \ll M_3$.
$\qed$

Recall that $P$ is a non-cyclic $p$-group where $p$ is a prime.
Assume that $\exp(P)>p$.
First, we define a subset
$\mathcal{S}'(P)$ of the set of all maximal cyclic subgroups of group
$P$ as follows:
\begin{align*}
C \in \mathcal{S}'(P)&\text{ if and only if }C \in \mathcal{M}(P),\ \lvert C\rvert >p\text{ and}\\
&\text{there is no }C'\in\mathcal M(P)\text{ of order greater than }\lvert C\rvert \text{ such that } C\ll C'.
\end{align*}
%
%
%
Note that for $C, C' \in \mathcal{S}'(P)$,
if $C \ll C'$, then we must have $\lvert C\rvert=\lvert C'\rvert$. As a result, the relation $\ll$ is symmetric on $\mathcal{S}'(P)$. By defining $\equiv$ as the reflexive closure of the domination relation on the set $\mathcal{S}'(P),$ we obtain an equivalence relation (by Lemma~\ref{ma-lem-1}, it is clear that  $\equiv$  is transitive).

Then, we define the set $\mathcal{S}(P)$ induced by $P$ as follows:

\begin{equation*}
\mathcal{S}(P)=\left\{
\begin{array}{ll}

\{\langle x \rangle\}, \hbox{where $x$ is any element of order $p$ in $P$,} & \hbox{if $\exp(P)=p$;} \\
\hbox{any system of representatives for the $\equiv$-classes}, & \hbox{if $\exp(P)>p$.}
\end{array}
                        \right.
\end{equation*}

We define $\mathcal{F}(P)$ as any subset of $P$ satisfying
\begin{equation}\label{subg-d-p}
\mathcal{F}(P)=\left\{
\begin{array}{ll}
\{x\}, & \hbox{if $\exp(P)=p$;} \\
\{x_1,x_2,\ldots,x_k\}, & \hbox{if $\exp(P)>p$,}
\end{array}
                        \right.
\end{equation}
where $x$ is an element of order $p$ of $P$, and
$\{\langle x_1\rangle,\langle x_2\rangle,\ldots,\langle x_k\rangle\}=\mathcal{S}(P)$. Note that $\mathcal{S}(P)$ and $\mathcal{F}(P)$ are actually families of sets; however, we will always consider particular representatives of these families.

We use the following example to illustrate the introduced sets.

\begin{example}\label{exam-1}
Consider $Q_{2^k}$ with $k\ge 3$, where $Q_{2^k}$ is the generalized quaternion group with $m=2^{k-2}$, as presented in \eqref{q4m}.
By \eqref{q4m-1} and \eqref{q4m-2}, we have that
\begin{equation*}
\mathcal{S}'(Q_{2^k})=\left\{
\begin{array}{ll}
\{\langle x\rangle,\langle xy\rangle,\langle x^2y\rangle\}, & \hbox{if $k=3$;} \\
\{\langle x\rangle\}, & \hbox{otherwise.}
\end{array}
                        \right.
\end{equation*}
As a result, for any $k\ge 3$, we have $\mathcal{S}(Q_{2^k})=\{\langle x\rangle\}$. Thus, we may choose $\mathcal{F}(Q_{2^k})=\{x\}$.
\end{example}

Now, consider again the non-cyclic nilpotent group $G$ which is not a $p$-group, as presented in \eqref{nil-fj}, and let $\mathbb{Z}_{n}=\langle x\rangle$. We naturally define
\begin{equation}\label{FG}
\mathcal{F}(G)=\left\{
\begin{array}{ll}
\bigcup_{i=1}^{t}\mathcal{F}(P_i), & \hbox{if $n=1$;} \\
\{x\}\cup (\bigcup_{i=1}^{t}\mathcal{F}(P_i)), & \hbox{otherwise.}
\end{array}
                        \right.
\end{equation}

\begin{lemma}\label{ni-kz-lem1}
Let $G$ be a non-cyclic nilpotent group which is not a $p$-group. Then $\mathcal{F}(G)$ is a dominating set of $\mathcal{D}(G)$.
\end{lemma}
\proof
Let $\mathcal{F}=\mathcal{F}(G)$ and $g\in V(\mathcal{D}(G))$. Let $G$ be presented as in \eqref{nil-fj}. Then, we may assume that
$$
g=g_1g_2\cdots g_tx^k,
$$
where $g_i \in P_i$ for each $1\le i \le t$, $x$ is a generator of $\mathbb Z_n$, and $k\in\mathbb N$. Note that for any $1\le i \le t$, there exists at least an element in $P_i\cap \mathcal{F}$.
Let $x_i\in P_i\cap \mathcal{F}$ for each $1\le i \le t$.
Note that $g$ and $x_i$ are non-trivial, where $1\le i \le t$.
If there exists an index $j\in \{1,2,\ldots,t\}$ such that $g_j=e$, then Lemma~\ref{two-con} implies that $x_j$ is adjacent to $g$ in $\mathcal{D}(G)$.

In the following, we may suppose that $g_i\ne e$ for each $1\le i \le t$. In view of Lemma~\ref{iso-vertex}, we see that $\langle g\rangle$ is not a maximal cyclic subgroup. Now, assume that $\langle g_i \rangle \in \mathcal{M}(P_i)$, for all $1\leq i\leq t$. In that case $\langle x^k \rangle \neq \langle x \rangle$, and it is easy to see that $g$ and $x$ are adjacent in $\mathcal{D}(G)$.

It remains to consider the case where some $g_i$ do not generate a maximal cyclic subgroup of $P_i$ and where none of $g_i$'s is the identity. Let $l\in\{1,2,\dots,t\}$ be such that $\langle g_l\rangle\not\in\mathcal M(P_l)$. Then, there is an element $y_l$ such that $\langle g_l\rangle< \langle y_l\rangle\in\mathcal M(P_l)$. Furthermore, there is some maximal cyclic subgroup of $P_l$ $\langle z_l\rangle\in\mathcal S(P_l)$ such that $\langle y_l\rangle\ll \langle z_l\rangle$. Since $\langle z_l\rangle$ contains all proper subgroups of $\langle y_l\rangle$, it follows that $\langle g_l\rangle< \langle z_l\rangle$. Finally, $g\sim z_l$. Without loss of generality, we may assume that $z_l\in\mathcal F(P_l)$. Thus, $\mathcal F$ is a dominating set of $\mathcal D(G)$.
$\qed$

\begin{thm}\label{nil-diam-bs}
Let $G$ be a finite nilpotent group which is not a $p$-group. Then $\mathcal{D}(G)$ is connected and $\diam(\mathcal{D}(G))\le 4$.
\end{thm}
\proof
Note that if $G$ is additionally non-cyclic, as presented in \eqref{nil-fj},  then for any two elements $a,b\in \mathcal{F}(G)$, we have $d(a,b)\leq 2$. Moreover, $d(a,b)=2$ if and only if $a,b\in P_i$, for some $1\leq i \leq t$. Therefore, in this case we have $\diam(\mathcal{F}(G))\leq 2$.

Now, based on Theorem~\ref{cyc-diam}, Lemma~\ref{ni-kz-lem1} and the conclusion above, the claim follows, which was also obtained in \cite{BC22} (see \cite[Theorem 5.3]{BC22}).
$\qed$

\begin{example}
Let us observe $Q_{2^k}\times \mathbb{Z}_n$, where $n\ge 3$ is an odd integer, $k\ge 3$, and $Q_{2^k}$ is the generalized quaternion group with $m=2^{k-2}$. Now by Theorem~\ref{d-n-1}, Example~\ref{exam-1} and Lemma~\ref{ni-kz-lem1}, we have $\gamma(\mathcal{D}(Q_{2^k}\times \mathbb{Z}_n))=2$.
\end{example}

\begin{lemma}\label{ni-kz-lem3}
Let $G$ be a non-cyclic nilpotent group that is not a $p$-group.
Then $\diam(\mathcal{D}(G))=4$ if and only if $G$ has two distinct cyclic subgroups $\langle u\rangle,\langle v\rangle \leq P$ that are not maximal cyclic in $P$, where $P$ is a Sylow $p$-subgroup of $G$, and $o(u)=o(v)$. Moreover, for  $a,b\in V(\mathcal{D}(G))$, we have $d(a,b)=4$ if and only if $a$ and $b$ have the following types:
$$a=uh,~~b=vh',$$
where $\langle h\rangle, \langle h'\rangle \in \mathcal{M}(H)$, and $G=P\times H.$
\end{lemma}
\proof
We first prove the sufficiency. Let $G$ be presented as in \eqref{nil-fj}.
Then, without loss of generality, we may assume that $\langle u\rangle,\langle v\rangle \leq P_1$, and $o(u)=o(v)$. Let $H$ be the product of all other Sylow subgroups of $G$, and let $h, h' \in H$  such that $\langle h \rangle, \langle h' \rangle \in \mathcal{M}(H)$.
Let us consider the elements
$$
a=uh, ~~b=vh'.
$$
Note that $a\ne u$ and $b\ne v$.
Since both $\langle u\rangle$ and $\langle v\rangle$ are not maximal cyclic in $G$, we may assume that $\langle u\rangle\lneq \langle u'\rangle\in \mathcal{M}(P_1)$ and $\langle v\rangle \lneq \langle v'\rangle\in \mathcal{M}(P_1)$. It is clear that $a$ and $b$ are non-adjacent. Also, in $\mathcal{D}(G)$, we have that
$$
N(a)=\{u_0\hat{h}| \langle u\rangle\lneq\langle u_0\rangle, o(u_0)=p_1^k,\langle \hat{h}\rangle \lneq \langle h\rangle,k\in \mathbb{N}\}
$$
and
$$
N(b)=\{v_0h''| \langle v\rangle\lneq\langle v_0\rangle, o(v_0)=p_1^l,\langle h''\rangle \lneq \langle h'\rangle,l\in \mathbb{N}\}.
$$
It follows that
$$
a\sim u' \sim h\sim v' \sim b.
$$
This means that
\begin{equation}\label{PM-eq-1}
d(a,b)\le 4.
\end{equation}

Since $\langle u \rangle$ and $\langle v \rangle$ are cyclic subgroups of the same order and $\langle u \rangle \neq \langle v \rangle$, they cannot both be contained in the same cyclic subgroup of group $P_1$.  Thus, $\langle u_0,v_0\rangle$ is not cyclic and $N(a)\cap N(b)=\emptyset$. It follows that $d(a,b)> 3$.

Now \eqref{PM-eq-1} implies that $d(a,b)=4$, and by Theorem~\ref{nil-diam-bs} we have that $\diam(\mathcal{D}(G))=4$, as desired.

We next prove the necessity. Suppose that $\diam(\mathcal{D}(G))=4$. Let $a$ and $b$ be two distinct vertices of $\mathcal{D}(G)$ with $d(a,b)=4$. Since we know that $\diam(\mathcal{F}(G))\leq 2$, the elements $a$ and $b$ cannot both belong to $\mathcal{F}(G)$. Thus, without loss of generality, we may assume that $a\in \mathcal{F}(G)$ and $b\notin \mathcal{F}(G)$.
In this case, we know that there exists an element $y\in \mathcal{F}(G)$ such that $y$ is adjacent to $b$. But now we have $d(a,b)\leq 3$, which leads to a contradiction.

Therefore, $a,b\notin \mathcal{F}(G)$. Combining $d(a,b)=4$ and Lemma~\ref{ni-kz-lem1}, we have that there exist distinct $z,w\in \mathcal{F}(G)$ such that $a\sim z$, $b\sim w$. Moreover, $z,w\in P_i$, for some $1\leq i\leq t$. Without loss of generality, we may assume that $i=1$.

For any $g\in V(\mathcal{D}(G))$, where $g\in P_j, 1<j\leq t$, we have that
\begin{equation}\label{PM-eq-2}
g\notin N(a)\cup N(b).
\end{equation}
For example, suppose that $g\in N(a)$. Then $a\sim g \sim w\sim b$, implying that $d(a,b)\leq 3$, which is a contradiction. A similar conclusion is obtained if $g\in N(b).$

Also, $p_1\mid o(a)$ and $p_1\mid o(b)$. If this were not the case, say $p_1\nmid o(a)$, then we would have $a\sim w\sim b$, which again leads to a contradiction.
Let $$\langle a\rangle=\langle u\rangle \times Q_1,~~\langle b\rangle=\langle v\rangle \times Q_2,$$
where $\langle u\rangle$ and $\langle v\rangle$ are the Sylow $p_1$-subgroups of $\langle a\rangle$ and $\langle b\rangle$, respectively. Since $a\sim z$ and $b\sim w$, it follows that $\langle u\rangle \lneq \langle z\rangle$ and
$\langle v\rangle \lneq \langle w\rangle$.
If $\langle u\rangle\leq\langle v\rangle$, then $a\sim w$, and so $d(a,b)\le 2$, which is a contradiction. As a result, $\langle u\rangle \nleq \langle v\rangle$.
Similarly, we also have $\langle v\rangle\nleq\langle u\rangle$. If $o(u)=o(v)$, then the required result holds.
Without loss of generality, assume that $o(u)<o(v)$. Then taking $\langle v'\rangle\leq \langle v\rangle$
with $o(u)=o(v')$, we have that $\langle u\rangle,\langle v'\rangle$ are two desired distinct nonmaximal cyclic subgroups of $P_1$.

Additionally, by \eqref{PM-eq-2}, it follows that $a$ and $b$ must be exactly of the form described in the formulation of the lemma itself.
$\qed$

\begin{lemma}\label{ni-kz-lem4}
Let $G$ be a non-cyclic nilpotent group which is not a $p$-group. Then $\diam(\mathcal{D}(G))=2$ if and only if $\exp(G)$ is a square-free positive integer.
\end{lemma}
\proof
We first prove the sufficiency. Denote $m=\exp(G)$, and suppose that $m$ is a square-free positive integer.
Then every non-trivial element of $G$ whose order is not equal to $m$ belongs to
$V(\mathcal{D}(G)).$

Let $a$ and $b$ be two distinct non-adjacent vertices of $\mathcal{D}(G)$. The existence of such elements is guaranteed by the fact that $G$ is a non-cyclic group with a square-free exponent. For example, there exist at least two elements of the same order that are non-adjacent in $\mathcal{D}(G)$. If
$\gcd(o(a),o(b))=1$, then $a$ and $b$ commute, and so they are adjacent, which contradicts the way these vertices were chosen.
Hence, $\gcd(o(a),o(b))=d$, where $1<d<m$.
Write
$$
o(a)=a_1\cdot d,~~o(b)=b_1\cdot d,
$$
where $a_1$ and $b_1$ are two positive integers. Since $\gcd(a_1,b_1)=1$, we can write $m=a_1 \cdot b_1 \cdot d \cdot c_1$, where $c_1$ is a positive integer. We observe that $a_1, b_1$ and $c_1$ cannot all be equal to 1.
Thus,
we may choose a non-trivial element $c\in G$ such that
$$
o(c)= a_1 \cdot b_1 \cdot  c_1, ~~|\langle c\rangle \cap \langle a \rangle|=a_1, ~~|\langle c\rangle \cap \langle b \rangle|=b_1.
$$
Both $\langle a,c\rangle$ and $\langle b,c\rangle$ are cyclic and $|\langle a,c\rangle|=|\langle b,c\rangle|=m$.
Note that $o(a)\ne m$ and $o(b)\ne m$.
As a result, there exist two prime divisors $p,q$ of $o(c)$ such that $p\nmid o(a)$ and $q\nmid o(b)$. Furthermore, note that $d \nmid o(c)$.
It follows that both $a$ and $b$ are adjacent to $c$ in $\mathcal{D}(G)$. Therefore, in this case, we have $d(a,b)=2$.
It follows that $\diam(\mathcal{D}(G))=2$, as desired.

We next prove the necessity. Suppose that $\diam(\mathcal{D}(G))=2$. Assume, to the contrary, that there exists a prime number $p$ such that $p^2\mid \exp(G)$. Let $g$ be a generator of a maximal cyclic Sylow $p$-subgroup of $G$, and denote by $H$ the product of all other Sylow subgroups of $G$.  Let now $h\in H$ be such that $\langle h\rangle \in \mathcal{M}(H)$.

Consider the elements $g^p$ and $g^ph$. It is easy to see that
$$
N(g^ph)=\{gh':  \langle h'\rangle \lneq \langle h \rangle \}.
$$

For every $h'$ such that $\langle h'\rangle \lneq \langle h \rangle$, we have  $\langle g^p\rangle \lneq \langle gh' \rangle$, so $g^p$ and $gh'$ are not adjacent. However, this implies that $d(g^p, g^ph)>2$, which is a contradiction.

$\qed$

Let $G$ be a non-cyclic nilpotent group which is not a $p$-group. Then $G$ is called a {\em $\Psi$-group} if it has two distinct cyclic subgroups $\langle u\rangle,\langle v\rangle \leq P$ that are not maximal cyclic in $P$, where $P$ is a Sylow $p$-subgroup of $G$, and $o(u)=o(v)$.

Now, combining Lemmas~\ref{ni-kz-lem3} and \ref{ni-kz-lem4}, and Corollary~\ref{nil-diam-bs}, we have the following theorem.

\begin{thm}\label{nil-g}
Let $G$ be a finite non-cyclic nilpotent group which is not a $p$-group. Then $\mathcal{D}(G)$ is connected. In particular,
\begin{equation*}
\diam(\mathcal{D}(G))=\left\{
                          \begin{array}{ll}
                          2, & \hbox{if $\exp(G) $ is a square-free number;} \\
                          4, & \hbox{if $G$ is a $\Psi$-group;}\\
                          3, & \hbox{otherwise.}
                          \end{array}
                        \right.
\end{equation*}

\end{thm}

Applying Theorem~\ref{nil-g} to abelian groups, we can obtain the  non-cyclic abelian groups whose difference
graph has diameter $4$. For example, by the following result, one can easily see that in \cite[Remark 5.4]{BC22}, $\diam(\mathcal{D}(\mathbb{Z}_4\times \mathbb{Z}_4 \times \mathbb{Z}_6))=4$.

\begin{cor}
Let $G$ be a finite non-cyclic abelian group which is not a $p$-group. Then
\begin{equation*}
\diam(\mathcal{D}(G))=\left\{
                          \begin{array}{ll}
                          2, & \hbox{if $\exp(G) $ is a square-free number;} \\
                          4, & \hbox{if $G$ has a subgroup isomorphic to  $\mathbb{Z}_{p^2}\times \mathbb{Z}_{p^2}$;}\\
                          3, & \hbox{otherwise,}
                          \end{array}
                        \right.
\end{equation*}
where $p$ is a prime.
\end{cor}

We now consider what can be said about the domination number of a non-cyclic nilpotent group. To this end, we recall the following elementary result.

\begin{lemma}
{\rm (\cite[Theorem~5.4.10~(ii)]{Gor})}\label{unp}
Given a prime $p$,
a $p$-group having a unique subgroup of order $p$ is either a cyclic $p$-group or a generalized quaternion $2$-group as presented in \eqref{q4m}.
\end{lemma}

Let $P$ be a $p$-group. Note that if $P$ is a generalized quaternion group, then $P$ has at least three cyclic subgroups of order $4$, each containing the unique involution of $P$. If $P$ is not a generalized quaternion group, then it contains at least three minimal subgroups. Accordingly, let $\mathcal{N}(P)$ denote the set of all cyclic subgroups of order $k$ in $P$, where $k$ is the smallest positive integer such that $P$ has more than one cyclic subgroup of order $k$.

\begin{lemma}
Let $G=P_1\times P_2 \times \mathbb{Z}_{n}$, where
$P_i$ is a nontrivial non-cyclic $p_i$-group for $i=1,2$, $p_1$ and $p_2$ are distinct primes, $n\ge 2$ and $p_i\nmid n$.
Then $\gamma(\mathcal{D}(G))\ge 3$.
\end{lemma}
\proof
By Corollary~\ref{nil-1}, we have $\gamma(\mathcal{D}(G))\ge 2$. Let $\mathbb{Z}_{n}=\langle g\rangle$.
It remains to show that $\gamma(\mathcal{D}(G))\ne 2$. Assume, for the sake of contradiction, that $\{a,b\}$ is a dominating set of $\mathcal{D}(G)$.

First, suppose that $\langle a\rangle$ and $\langle b\rangle$ both have subgroups belonging to $\mathcal{N}(P_1)$. Then there exists an element $u$ such that $\langle u\rangle \in \mathcal{N}(P_1)$ and $\langle u\rangle \nleq \langle a\rangle, \langle b\rangle$. However, in this case $u$ is adjacent to neither $a$ nor $b$, contradicting the assumption that $\{a,b\}$ is a dominating set.

Now suppose, without loss of generality, that $\langle a\rangle$ contains a subgroup in $\mathcal{N}(P_1)$ while $\langle b\rangle$ does not, that is, $\Sub(\langle a\rangle)\cap \mathcal{N}(P_1)\neq \emptyset$ and $\Sub(\langle b\rangle)\cap \mathcal{N}(P_1)= \emptyset$. If $\langle b\rangle$ is not a maximal cyclic subgroup of $P_2\times \mathbb{Z}_n$, then we can choose $u$ such that $\langle u \rangle \in \mathcal{N}(P_1)$ and $\langle u\rangle \nleq \langle a\rangle$. In this case, the element $ub$ is adjacent to neither $a$ nor $b$, resulting in a contradiction. Therefore, $\langle b\rangle \in \mathcal{M}(P_2\times \mathbb{Z}_n)$, and $\langle b\rangle$ contains a subgroup from $\mathcal{N}(P_2)$. Now take $u\in P_1$ such that $\langle u\rangle \in \mathcal{N}(P_1)$ and $u\notin \langle a\rangle$, and $v\in P_2$ such that $\langle v\rangle \in \mathcal{N}(P_2)$ and $v\notin \langle b\rangle$. Then $uv$ is adjacent to neither $a$ nor $b$, again contradicting our assumption.

It follows that $\Sub(\langle a\rangle)\cap \mathcal{N}(P_1)= \emptyset$ and $\Sub(\langle b\rangle)\cap \mathcal{N}(P_1)= \emptyset$. Similarly, $\Sub(\langle a\rangle)\cap \mathcal{N}(P_2)= \emptyset$ and $\Sub(\langle b\rangle)\cap \mathcal{N}(P_2)= \emptyset$. Thus, if neither $P_1$ nor $P_2$ is a generalized quaternion group, then $a,b\in \langle g\rangle$. If, without loss of generality, $P_1$ is a generalized quaternion group with unique involution $h$, then $a,b\in \langle g,h\rangle$. Let $w\in P_1$ be such that $\langle w\rangle \in \mathcal{N}(P_1)$. In both cases, the element $wg$ is adjacent to neither $a$ nor $b$, again contradicting our assumption. Therefore, $\{a,b\}$ cannot be a dominating set of $\mathcal{D}(G)$.
$\qed$





\section*{Acknowledgements}

The first author was supported by National Natural Science Foundation of China (Grants No. 11801441 \& 12326333) and  Shaanxi Fundamental Science Research Project for Mathematics and Physics (Grant No. 22JSQ024). 
The second author acknowledges the financial support of the
Ministry of Science, Technological Development and Innovation of the
Republic of Serbia (Grants No. 451-03-137/2025-03/200125 \&
451-03-136/2025-03/200125). The third author was supported by the Serbian Ministry of Science, Technological Development and Innovation through the Mathematical Institute of the Serbian Academy of Sciences and Arts.

\end{document}